\documentclass[tbtags,a4paper,11pt]{amsart}
\usepackage{xspace}
\usepackage{amsmath}
\usepackage{amssymb,a4}

\newcommand\R{\ensuremath{\mathbb{R}}\xspace}  
\newcommand\N{\ensuremath{\mathbb{N}}\xspace}
\newcommand{\aop}{\operatorname{\mathbf a}}  
\newcommand{\diw}{\ensuremath{\mathrm{div}}}
\newcommand{\capp}{\ensuremath{\text{cap}_p}}

\newcommand\dx{\ensuremath{\dleb\!x}}
\newcommand\dleb{\operatorname{d}}
\newcommand\mlimsup{\operatornamewithlimits{\overline{lim}}}
\newcommand\indi{\nbOne}
\newcommand\nbOne{{\mathchoice {\mathrm {1\mskip-4.1mu l}} {\mathrm{ 1\mskip-4.1mu
        l}} {\mathrm {1\mskip-4.6mu l}} {\mathrm {1 \mskip-5.2mu l}}}}
\newcommand\psidp{\psi_\delta^+}
\newcommand\psidm{\psi_\delta^-}
\newcommand\OU{{\mathcal U}}

\author[{O. Guib\smash{\'e}}]
{Olivier Guib\'e}

\title[Remarks on the
  uniqueness of comparable solutions \dots]{Remarks on the
  uniqueness of comparable renormalized solutions of
elliptic equations with measure data}

\makeatletter
\def\th@definition{%
  \normalfont 
}
\def\th@plain{%
  \slshape 
}
\def\th@remark{%
  \normalfont 
  \thm@preskip\topsep
  \divide\thm@preskip\tw@
  \thm@postskip\thm@preskip
}
\makeatother
\theoremstyle{remark}
\newtheorem{rmk}{Remark}

\theoremstyle{definition}
\newtheorem{defi}[rmk]{Definition}
\theoremstyle{plain}
\newtheorem{prop}[rmk]{Proposition}

\newtheorem{theorem}[rmk]{Theorem}

\begin{document}

\maketitle
\ \\ \vglue-1.5cm
\begin{minipage}{0.9\linewidth}
\begin{center}
{Laboratoire de~Math\'ematiques~Rapha\"el~Salem}\\ 
UMR CNRS 6085, Site Colbert \\
Universit\'e de Rouen\\
F-76821 Mont Saint Aignan cedex \\
 E-mail : \texttt{olivier.guibe@univ-rouen.fr} 
\end{center}
\end{minipage}

\thispagestyle{empty}
\ \\
\vglue 2cm 

\begin{abstract}
  We give a partial uniqueness result concerning comparable
  renormalized solutions of the nonlinear elliptic problem
  $-\diw(\aop(x,Du))=\mu$ in $\Omega$, $u=0$ on $\partial\Omega$,
  where $\mu$ is a Radon measure with bounded variation on $\Omega$.
\end{abstract}

\renewcommand\baselinestretch{1.1}
\normalsize

\medskip
\medskip

\setcounter{equation}{0}
\section{Introduction}
Let us consider the nonlinear elliptic problem
\begin{align}
  -\diw\big(\aop(x,Du)\big)= \mu
  & \quad\text{ in } \Omega, \label{og1}\\
u= 0 &\quad \text{ on } \partial\Omega,\label{og2}
\end{align}
where $\Omega$ is a bounded open subset of $\R^N$ with $N\geq 2$,
$u\mapsto -\diw\big(\aop(x,Du)\big)$ is a strictly monotone operator
from $W^{1,p}_0(\Omega)$ into $W^{-1,p'}(\Omega)$ and $\mu$ is a Radon
measure with bounded variation on $\Omega$.
\par
In the linear case G. Stampacchia has defined in \cite{St} the notion of
``solution by transposition'' which insures existence and uniqueness
of such a solution. If $p=2$ and for the nonlinear case, this notion
is generalized in \cite{M94} and the existence and uniqueness of the solution
obtained as limit of approximations is proved in \cite{M94} (see also
\cite{Bo97} and, for a class of pseudo-monotone operator \cite{DG}).
\par
If $2-1/N<p\leq N$ the existence of a solution of
(\ref{og1})--(\ref{og2}) in the sense of distributions is proved by
L.~Boccardo and T.~Gallou\"et in \cite{BG89}. However,  using
the counter example of J.~Serrin \cite{serrin} it is well known that this solution is
not unique in general, except in the case $p=N$ for an appropriate
choice of the space to which the solution belongs (see \cite{DHM99}
and  \cite{GIS97}). 
\par
When $\mu$ is a function of $L^1(\Omega)$  the notions of entropy
solution \cite{BBGGPV95}, of solution obtained as limit of approximations
\cite{Dall96} and of renormalized solution \cite{LM} (see also
\cite{M93} and \cite{M94}) provide existence and
uniqueness results (and these three notions are actually equivalent).
\par
When $\mu$ is a Radon measure with bounded variation on $\Omega$,
G.~Dal~Maso, F.~Murat, L.~Orsina and A.~Prignet have recently introduced
 in \cite{DMOP} and \cite{DMOP99} a notion of renormalized
solution of (\ref{og1})--(\ref{og2}) which generalizes the three
(and equivalent) previous ones. The authors prove in \cite{DMOP99} the existence of
such a renormalized solution, a stability result and partial
uniqueness results for ``comparable'' solutions.
In particular, under some assumptions on $\aop$,  if $u_1$ and $u_2$
are two renormalized solutions of (\ref{og1})--(\ref{og2}) such that
$u_1-u_2$ belongs to $L^\infty(\Omega)$ (this condition is here the
precise meaning of the fact that the two solutions are comparable), then $u_1=u_2$. 
 The uniqueness of the renormalized solution
of (\ref{og1})--(\ref{og2}) remains an open problem in general
 and the present paper is devoted to weaken this 
condition. We prove that the condition of being comparable can be localized in a
neighborhood $\OU$ of the set where $\mu$ is singular and that it is sufficient to
assume that $(u_1-u_2)^-$ (the negative part of $u_1-u_2$) belongs to $L^\infty(\OU)$.
\par\smallskip
The paper is organized as follows: Section 2 is devoted to give the
assumptions on the data and to recall the definition of a renormalized
solution of (\ref{og1})--(\ref{og2}). In Section 3 (Theorems
\ref{ogt1} and \ref{ogt2}) we establish
partial uniqueness results concerning comparable renormalized
solutions of (\ref{og1})--(\ref{og2}).

\section{Assumptions and definitions}
Let $\Omega$ be a bounded open subset of $\R^N$ with $N\geq 2$, $p$
and $p'$ two real numbers such that $1<p<N$ and $1/p+1/{p'}=1$. We
assume that $\aop\,:\,\Omega\times\R^N\mapsto \R^N$ is a
Carath\'eodory function (i.e. measurable with respect to $x$ and
continuous with respect to $\xi$) such that
\begin{equation}
  \aop(x,\xi)\cdot \xi \geq \alpha |\xi|^p, \label{h1a} 
\end{equation}
\begin{equation}
  \big(\aop(x,\xi)-\aop(x,\xi')\big)\cdot(\xi-\xi') > 0,\label{h1b}
\end{equation}
\begin{equation}
  \label{h2}
\big|\aop(x,\xi)\big|\leq \gamma\big( b(x)+ |\xi|^{p-1}\big)
\end{equation}
for every $\xi$, $\xi'$ ($\xi\neq\xi'$) in $\R^N$ and almost everywhere
in $\Omega$, where $\gamma>0$, $\alpha>0$ and $b$ is a nonnegative
function lying in $L^p(\Omega)$.
\par
We denote by ${\mathcal M}_b(\Omega)$ the set of Radon measures on
$\Omega$ with total bounded variation on $\Omega$ and by ${\mathcal M}_0(\Omega)$
the set of measures of ${\mathcal M}_b(\Omega)$ that are absolutely
continuous with respect to the $p$--capacity (i.e. $\mu\in{\mathcal
  M}_b(\Omega)$ and $\mu(E)=0$ for every Borel set $E$ such that
$\capp(E,\Omega)=0$). For $K>0$ we define as
$T_K(r)=\max(-K,\min(K,r))$ the truncation function at height $\pm K$.
If $A$ is a measurable set we denote by $\indi_A$ the characteristic
function of $A$. 
\par
We recall now a decomposition result of the  Radon measures (see
\cite{BGO} and \cite{FST}) and the definition of the gradient of a
function whose truncates belong to $W^{1,p}_0(\Omega)$ (see
\cite{BBGGPV95} Lemma 2.1 and \cite{LM}) which are needed to define  (following
\cite{DMOP99}) a renormalized solution of (\ref{og1})--(\ref{og2}).
\begin{prop}{\upshape (\cite{BGO} and \cite{FST})}\label{ogp1}
  Let $\mu$ be an element of ${\mathcal M}_b(\Omega)$. There exists
 two functions $f\in L^1(\Omega)$, $g\in(L^{p'}(\Omega))^N$, two nonnegative
  measures in ${\mathcal M}_b(\Omega)$, $\lambda^+$ and $\lambda^-$,
  which are concentrated respectively on two disjoint Borel sets $E^+$
  and $E^-$ of zero $p$-capacity such that
  \begin{equation*}
    \mu= f-\diw(g) +\lambda^+-\lambda^-.
  \end{equation*}
Moreover, if $\mu_0$ denotes $f-\diw(g)$ then $\mu_0\in {\mathcal
  M}_b(\Omega)$ and the decomposition $\mu=\mu_0+\lambda^+-\lambda^-$
  is unique.
\end{prop}
\par\smallskip
\begin{defi}{ (\cite{BBGGPV95} and \cite{LM})}\label{ogd1}
  Let $u$ be a measurable function defined from $\Omega$ into
  $\overline{\R}$ which is
  finite almost everywhere in $\Omega$. Assume that $T_K(u)\in
  W^{1,p}_0(\Omega)$ $\forall K>0$. Then there exists a unique
  measurable function 
  $v\,:\,\Omega\mapsto \R^N$ such that
  \begin{equation*}
    \forall K>0,\quad DT_K(u)=\indi_{\{|u|<K\}} v \quad\text{a.e. in }\Omega.
  \end{equation*}
This function $v$ is called the gradient of $u$ and is denoted by $Du$.
\end{defi}
\par\smallskip

Following \cite{DMOP99} we are now in a position to recall the definition of renormalized
solution.
\begin{defi}(\cite{DMOP99})\label{ogd2}
  Let $\mu$ be an element of ${\mathcal M}_b(\Omega)$ and
  $\mu=f-\diw(g)+\lambda^+-\lambda^-$ the decomposition given by
  Proposition \ref{ogp1}.
A function $u$ defined from $\Omega$ into $\overline{\R}$ is a renormalized
solution of (\ref{og1})--(\ref{og2}) if
 \begin{equation}\label{ogd2n1}
   \begin{minipage}[t]{0.9\linewidth}
     $u$ is measurable and finite almost everywhere in $\Omega$ and
     $T_K(u)\in W^{1,p}_0(\Omega)$ $\forall K>0$;
   \end{minipage}
 \end{equation}
 \begin{equation}
   \label{ogd2n2}
   |Du|^{p-1}\in L^q(\Omega)\quad \forall q<\frac{N}{N-1};
 \end{equation}
 \begin{gather}
   \label{ogd2n3}
   \begin{minipage}[t]{0.9\linewidth}
$\forall w\in W^{1,p}_0(\Omega)\cap L^\infty(\Omega)$ such that
$\exists K>0$ and two functions $w^{+\infty}$ and $w^{-\infty}$ lying in
$W^{1,r}(\Omega)\cap L^\infty(\Omega)$ with $r>N$ and 
\end{minipage} \\
\left\{\begin{aligned}
  w=w^{+\infty} & \text{ on $\{ x\,;\, u(x)>K\}$}, \\
 w=w^{-\infty} & \text{ on $\{x\,;\, u(x)<-K\}$}, 
\end{aligned}\right. \notag 
\end{gather}
we have
\begin{multline}  \int_\Omega \aop(x,Du)\cdot Dw \dx =   \int_\Omega f w \dx +
  \int_\Omega g\cdot Dw \dx  \\
  + \int_\Omega w^{+\infty} \dleb\!\lambda^+ -
  \int_\Omega w^{-\infty} \dleb\! \lambda^-. 
\label{ogd2n4}
 \end{multline}
\end{defi}

\par\medskip
It is proved in \cite{DMOP99} that if $\aop$ verifies (\ref{h1a}),
(\ref{h1b}) and (\ref{h2}) then for any element $\mu$ belonging to
${\mathcal M}_b(\Omega)$ there exists at least a renormalized solution
of (\ref{og1})--(\ref{og2}). 
\par\smallskip
\begin{rmk}\label{ogrmk1}
Every function $w\in {\mathcal C}^\infty_c(\Omega)$ is an admissible test
function in (\ref{ogd2n4}) and then any renormalized solution of
(\ref{og1})--(\ref{og2}) is also solution in the sense of
distributions.
\par
Furthermore if $\varphi\in W^{1,r}(\Omega)\cap L^\infty(\Omega)$ with
$r>N$ then we have
\begin{equation}
  \label{ogrn1}
  \lim_{n\rightarrow +\infty} \frac{1}{n} \int_{\{|u|<n\}}
  \aop(x,Du)\cdot Du \varphi \dx = \int_\Omega\, \varphi  \dleb\!
  \lambda^+
+ \int_\Omega \varphi  \dleb\!\lambda^-.
\end{equation}
  This property (see \cite{DMOP99} for more details on the properties
of renormalized solutions) is obtained by using the admissible test
function $w=\frac{1}{n} T_n(u) \varphi$ in (\ref{ogd2n4}) and by passing
to the limit as $n$ goes to infinity.
\end{rmk}

\section{Uniqueness of comparable solutions}

In \cite{DMOP99} the authors prove under assumptions (\ref{h1a})
and (\ref{h1b}),  the strong monotonicity of $\aop$ and the local
Lipschitz continuity, or the H\"older continuity, with respect to
$\xi$, i.e. $\aop$ verifies
\begin{equation}
  \label{ogt2h}
  \left\{ 
    \begin{aligned}
      \big(\aop(x,\xi)-\aop(x,\xi')\big)\cdot (\xi-\xi') & \geq \alpha
      |\xi-\xi'|^p & \quad \text{ if } p\geq 2 \\[0.2cm]
 \big(\aop(x,\xi)-\aop(x,\xi')\big)\cdot (\xi-\xi') & \geq \alpha
     \smash{ \frac { |\xi-\xi'|^2}{\big(|\xi|+|\xi'|\big)^{2-p}}} & \quad \text{ if
      } p<2, 
    \end{aligned}\right.
\end{equation}
\begin{equation}
  \label{ogt1h0}
\left\{\begin{aligned}
      \big|\aop(x,\xi)-\aop(x,\xi')\big| & \leq \gamma
     \big(b(x)+|\xi|+|\xi'|\big)^{p-2} |\xi-\xi'| & & \quad \text{ if } p\geq 2, \\[0.1cm]
 \big|\aop(x,\xi)-\aop(x,\xi')\big| & \leq \gamma
   |\xi-\xi'|^{p-1}  & & \quad \text{ if
      } p<2,
    \end{aligned}\right.
\end{equation}
for every $\xi$, $\xi'\in\R^N$ and almost everywhere in $\Omega$,
where $\gamma>0$ and $b$ is a nonnegative function in $L^p(\Omega)$, 
that if two renormalized solutions $u_1$ and $u_2$ of (\ref{og1})--(\ref{og2})
(relative to the same element $\mu \in {\mathcal M}_b(\Omega)$)  satisfy
the condition of being comparable, in the sense that $u_1-u_2\in
L^\infty(\Omega)$, then 
$u_1=u_2$. In Theorem~\ref{ogt2} below we weaken this condition; if
there exists an open neighborhood $\OU$ of $E=E^+\cup E^-$ where
$E^+$ and $E^-$ are given by Proposition \ref{ogp1} such that
$(u_1-u_2)^-\in L^\infty(\OU)$, then $u_1=u_2$. This result is  a
consequence of the following theorem.

\begin{theorem}\label{ogt1}
Assume that (\ref{h1a}), (\ref{h1b}), (\ref{h2}) and (\ref{ogt1h0}) hold true. Let
$\mu$ be an element of ${\mathcal M}_b(\Omega)$ and let $E=E^+\cup E^-$
where $E^+$ and $E^-$ are the two disjoint Borel sets of zero
$p$-capacity given by Proposition \ref{ogp1}. Let $u_1$ and $u_2$ be
two renormalized solutions of (\ref{og1})--(\ref{og2}) with $\mu$ as
right-hand side. If moreover there exists an open set $\OU$ such that
\begin{gather}
  \label{ogt1h1}
  E\subset \OU \subset \Omega, \\
\forall K>0\quad \lim_{n\rightarrow +\infty} \frac{1}{n} \int_{\substack{
{\mathcal U}\cap \{u_1-u_2<K\} \\
\cap \{|u_1|<n,\, |u_2|<n\}
}}
\big|Du_1-Du_2|^p \dx =0, \label{ogt1h2}
\end{gather}
then $u_1=u_2$.
\end{theorem}

\begin{rmk}
  Using the following property for every $m\in\N^*$
  \begin{equation*}
    \begin{split}
      \{|u_1|<2^m,\, |u_2|<2^m\} \subset & \{|u_1|<1,\,|u_2|<1\}  \\
& {} \cup \bigcup_{k=0}^{m-1}   \{2^k\leq |u_1|<2^{k+1},\,|u_2|<2^{k+1}\}
\\
& {} \cup \bigcup_{k=0}^{m-1}   \{2^k\leq |u_2|<2^{k+1},\,|u_1|<2^{k+1}\},
\end{split}
  \end{equation*}
 a Cesaro  argument and  the
  fact that $T_1(u_1)$ and $T_1(u_2)$ belong to $W^{1,p}_0(\Omega)$,
  the condition (\ref{ogt1h2})  is
  equivalent to
\begin{multline*}
  \lim_{n\rightarrow +\infty}\frac{1}{n}\left( \int_{\substack{
{\mathcal U}\cap \{u_1-u_2<K\} \\
\cap \{n\leq |u_1|<2n,\, |u_2|<2n\} 
}}
\big|Du_1-Du_2|^p \dx  \right. \\
\left. + \int_{\substack{
{\mathcal U}\cap \{u_1-u_2<K\} \\
\cap \{|u_1|<2n,\, n\leq |u_2|<2n\} 
}}
\big|Du_1-Du_2|^p \dx \right)= 0,
\end{multline*}
for all $K>0$. Notice that the condition above with $\OU=\Omega$ and
$K=+\infty$ (so that $\OU\cap \{u_1-u_2<K\}=\Omega$) is the one given in
\cite{DMOP99} (Theorem 10.3).
\end{rmk}
\par
\medskip
\begin{theorem}\label{ogt2}
  Assume that (\ref{h1a}), (\ref{h1b}), (\ref{h2}), (\ref{ogt2h}) and (\ref{ogt1h0})
  hold true. Let $\mu$ be an element  of ${\mathcal M}_b(\Omega)$ and let $E=E^+\cup E^-$
where $E^+$ and $E^-$ are the two disjoint Borel sets of zero
$p$-capacity given by Proposition \ref{ogp1}. Let $u_1$ and $u_2$ be
two renormalized solutions of (\ref{og1})--(\ref{og2}) with $\mu$ as
right-hand side. If moreover there exists an open set $\OU$ such that
  $E\subset\OU$ and $(u_1-u_2)^-\in L^\infty(\OU)$ (or
  $(u_1-u_2)^+\in L^\infty(\OU)$), then $u_1=u_2$.
\end{theorem}
\par
\medskip

\begin{proof}[Proof of Theorem \ref{ogt1}]
Using Proposition \ref{ogp1}, let $f\in L^1(\Omega)$, $g\in
(L^{p'}(\Omega))^N$, $\lambda^+$ and $\lambda^-$ two nonnegative
measures of ${\mathcal M}_b(\Omega)$ which are concentrated on two
disjoint subsets $E^+$ and $E^-$ of zero $p$-capacity such that $\mu=
f -\diw (g) +\lambda^+-\lambda^-$. Since $\capp(E^+,\Omega)=0$ and
$E^+\subset \OU \subset \Omega$ we have (see \cite{KKK})
$\capp(E^+,\OU)=0$ (and also $\capp(E^-,\OU)=0$). Thus, following the
construction of the cut-off functions in \cite{DMOP99}, we define for all
$\delta>0$ two functions, $\psidp$ and $\psidm$, lying in 
${\mathcal C}^\infty_c(\OU)$ such that
\begin{gather}
  0\leq \psidp \leq 1,\quad 0\leq \psidm\leq 1\quad\text{ on
  $\OU$,} \label{ogt1n1} \\
\text{supp}(\psidp)\cap \text{supp}(\psidm) =\emptyset, \label{ogt1n2} \\
\int_\Omega \psidm \dleb\!\lambda^+ < \delta,\quad 
\int_\Omega \psidp \dleb\!\lambda^- < \delta, \label{ogt1n3} \\
\int_\Omega (1- \psidp) \dleb\!\lambda^+< \delta, \quad 
\int_\Omega (1- \psidm) \dleb\!\lambda^-< \delta. \label{ogt1n4}
\end{gather}
Since $\OU\subset\Omega$, we define $\psidp\equiv \psidm\equiv 0$ on
$\Omega\setminus\OU$ so that we have $\psidp$, $\psidm\in {\mathcal
  C}^\infty_c(\Omega)$.
\par
For any $n\in\N^*$ let $h_n$ be the function defined by
$h_n(r)=\big(n-T_n^+(|r|-n)\big)/n$ $\forall r\in\R$.
\par
Let $K>0$ be fixed, $n\in\N^*$ and $\delta>0$. Since the function
$h_n$ belongs to $W^{1,\infty}(\R)$ while $\text{supp}(h_n)=[-2n,2n]$ is compact,
 from the regularity of $u_1$ and $u_2$ we obtain that 
the function $h_n(u_1)h_n(u_2) \big(T_K(u_1-u_2)-K(\psidp+\psidm)\big)$ lies in
$W^{1,p}_0(\Omega)\cap L^\infty(\Omega)$ and is equal to zero on the
set $\{x\,;\,|u_i(x)|>2n\}$ for $i=1,2$. Therefore setting
$W_K=T_K(u_1-u_2)$ the function $h_n(u_1)h_n(u_2)
\big(W_K-K(\psidp+\psidm)\big)$ is an admissible test 
function on both equations (\ref{og1}) written for $u_1$ and $u_2$,
relative to (\ref{ogd2n4}) of Definition \ref{ogd2}. Subtracting the resulting
equalities gives
\begin{align}
&\int_\Omega h_n(u_1) h_n(u_2) 
  \big(\aop(Du_1)-\aop(Du_2)\big)\cdot 
(DW_K-KD(\psidp-\psidm)) \dx
  \tag{{A}} \label{ogt1n5} \\
+&\int_\Omega h'_n(u_1)h_n(u_2)
\big(\aop(Du_1)-\aop(Du_2)\big)\cdot Du_1
\big(W_K-K(\psidp+\psidm)\big) \dx \tag{{B}}\label{ogt1n6} \\
+&\int_\Omega  h'_n(u_2)h_n(u_1)
\big(\aop(Du_1)-\aop(Du_2)\big)\cdot Du_2
\big(W_K-K(\psidp+\psidm)\big)  \dx  \tag{{C}} \label{ogt1n7} \\
= &\ 0 \notag
\end{align}
\par
In order to study the behavior of the terms above as $n$ goes to
infinity and $\delta$ goes to zero, \ref{ogt1n5} and \ref{ogt1n6}
are split into $\mathrm{A}_1+\mathrm{A}_2$ and
$\mathrm{B}_1+\mathrm{B}_2$ respectively, where
\begin{align*}\
  \mathrm{A_1} & =\int_\Omega h_n(u_1) h_n(u_2)
  \big(\aop(x,Du_1)-\aop(x,Du_2)\big)\cdot DT_K(u_1-u_2) \dx, \\
\mathrm{A_2} & = -K \int_\Omega h_n(u_1) h_n(u_2)
  \big(\aop(x,Du_1)-\aop(x,Du_2)\big) \cdot (D\psidp+D\psidm) \dx, \\
\mathrm{B_1} & = \int_\Omega h'_n(u_1) h_n(u_2)  \big(\aop(x,Du_1)-\aop(x,Du_2)\big)
  \cdot Du_1 W_K (1-\psidp-\psidm) \dx, \\
\mathrm{B_2} & = \int_\Omega h'_n(u_1) h_n(u_2)  \big(\aop(x,Du_1)-\aop(x,Du_2)\big)
  \cdot Du_1 \big( W_K-K\big) (\psidp+\psidm) \dx.
\end{align*}
\par

 From (\ref{h2}) and (\ref{ogd2n2}) it follows that $\aop(x,Du_i)$
belongs in particular to $L^1(\Omega)$  for $i=1,2$ and then
$\big(\aop(x,Du_1)-\aop(x,Du_2)\big)\cdot (D\psidp+D\psidm)$ belongs
to $L^1(\Omega)$. Since $h_n(u_1)h_n(u_2)$ converges to 1 almost
everywhere as $n$ tends to infinity and is uniformly bounded, Lebesgue
Theorem leads to 
\begin{equation*}
  \lim_{n\rightarrow+\infty} \mathrm{A_2}= -K \int_\Omega
  \big(\aop(x,Du_1)-\aop(x,Du_2)\big) \cdot (D\psidp+D\psidm)  \dx.
\end{equation*}
Recalling that $u_1$ and $u_2$ are also solution of
(\ref{og1})--(\ref{og2}) in the sense of distributions and since
$\psidp$, $\psidm\in {\mathcal C}^\infty_c(\Omega)$,  we obtain that
\begin{equation}\label{ogt1n8}
   \lim_{n\rightarrow+\infty} \mathrm{A_2}=0.
\end{equation}
\par
Due to the definition of $\psidp$ and $\psidm$ we have $1\geq
1-\psidp-\psidm\geq 0$. Thus Assumption (\ref{h2}) and Young's inequality
lead to
\begin{multline*}
  |\mathrm{B_1}| \leq \frac{C}{n}\left( 
\int_{\{|u_1|<2n\}} |Du_1|^p (
  1-\psidp-\psidm) \dx \right. \\
\left. {}+ \int_{\{|u_2|<2n\}}  |Du_2|^p (
  1-\psidp-\psidm) \dx
+\int_\Omega b^p \dx \right)
\end{multline*}
and (\ref{h1a}) gives
\begin{multline}\label{ogt1n9}
  |\mathrm{B_1}| \leq \frac{C}{n} \left(\int_{\{|u_1|<2n\}}
  \aop(x,Du_1)\cdot Du_1 (
  1-\psidp-\psidm) \dx \right. \\
\left. + \int_{\{|u_2|<2n\}}  \aop(x,Du_2)\cdot Du_2 (
  1-\psidp-\psidm) \dx
+\int_\Omega b^p \dx
\right),
\end{multline}
where $C$ is a generic constant independent of $n$ and $\delta$. Since
$1-\psidp-\psidm\in {\mathcal C}^\infty(\overline{\Omega})$, 
using the property (\ref{ogrn1}) of renormalized solutions we get
for $i=1,2$
\begin{multline*}
  \lim_{n\rightarrow+\infty} \frac{1}{2n} \int_{\{|u_i|<2n\}} \aop(x,Du_i)\cdot
  Du_i (1-\psidp-\psidm)\dx = \\
 \int_\Omega (1-\psidp-\psidm) \dleb\!\lambda^+
  + \int_\Omega (1-\psidp-\psidm) \dleb\!\lambda^-,
\end{multline*}
from which it follows, using (\ref{ogt1n3}), (\ref{ogt1n4}) and
(\ref{ogt1n9}) and since $b\in L^p(\Omega)$,
\begin{equation*}
  \mlimsup_{n\rightarrow+\infty}|\mathrm{B_1}|  \leq C \delta,
\end{equation*}
and then
\begin{equation}\label{ogt1n10}
  \lim_{\delta\rightarrow 0} \mlimsup_{n\rightarrow +\infty} |\mathrm{B_1}| =0.
\end{equation}
\par
Let  ${\mathcal U}_{n,K}$ be the set defined by
\begin{equation}\label{ogt1n11}
 {\mathcal  U}_{n,K}=  \OU\cap \{ |u_1|<2n\}\cap  \{
  |u_2|<2n\}  \cap \{ u_1-u_2<K\}.
\end{equation}
Because  $0\leq K-T_K(u_1-u_2)\leq 2 K \indi_{\{u_1-u_2<K\}}$ (recall
that $W_K=T_K(u_1-u_2)$), from the definition of the cut-off functions
$\psidp$ and $\psidm$ we obtain 
\begin{equation*}
    |\mathrm{B_2}| \leq \frac{2K}{n} \int_{{\mathcal  U}_{n,K}}
    \big|\aop(x,Du_1)-\aop(x,Du_2)\big| \,|Du_1| \dx.
\end{equation*}
Using H\"older inequalities together with (\ref{h2}) permits us to
deduce that if $p\geq 2$ then
\begin{multline*}
   |\mathrm{B_2}| \leq C K \left(\frac{1}{n} \int_{{\mathcal  U}_{n,K}}
   |Du_1-Du_2|^p \dx \right)^{1/p} \\
{} \times\left( \frac{1}{n}
\int_{{ \substack{\{|u_1|<2n \}\\\!\cap\{ |u_2|<2n \}}}}
   \big( b(x)+|Du_1|+|Du_2|\big)^p \dx \right)^{1/{p'}}
\end{multline*}
and if $p<2$ then
\begin{gather*}
   |\mathrm{B_2}| \leq {C K} \left(\frac{1}{n}\int_{{\mathcal
   U}_{n,K}}
    |Du_1-Du_2|^p \dx \right)^{1/{p'}} \left(\frac{1}{n}
\int_{{ \substack{\{|u_1|<2n \}\\\!\cap\{ |u_2|<2n \}}}}
 |Du_1|^p \dx \right)^{1/p}.
\end{gather*}
In both cases, property (\ref{ogrn1}) (with $\varphi\equiv 1$) and
(\ref{ogt1h2}) lead to
\begin{equation*}
  \forall \delta>0 \quad  \lim_{n\rightarrow +\infty} |\mathrm{B}_2| = 0.
\end{equation*}
From (\ref{ogt1n10}) it follows that 
\begin{gather}\label{ogt1n12}
  \lim_{\delta\rightarrow 0} \mlimsup_{n\rightarrow +\infty} |\mathrm{B}|=0
\quad\text{and}\quad
  \lim_{\delta\rightarrow 0} \mlimsup_{n\rightarrow +\infty}
  |\mathrm{C}|=0 \quad\text{(by symmetry).}
\end{gather}
\par
From  (\ref{ogt1n8}) and (\ref{ogt1n12}) we then have
 $\lim_{\delta\rightarrow 0} \lim_{n\rightarrow+\infty}
\mathrm{A_1}=0$. Since $h_n(u_1)h_n(u_2)$ is nonnegative and converges
to 1 almost everywhere in $\Omega$, the monotone character of the
operator $\aop$ and Fatou lemma imply that for all $K>0$
\begin{equation*}
  \int_{\{|u_1-u_2|<K\}} \big(\aop(x,Du_1)-\aop(x,Du_2)\big)\cdot
  (Du_1-Du_2) \dx = 0,
\end{equation*}
and from (\ref{h1b}) we can conclude that $u_1=u_2$.
\end{proof}

\par\smallskip

\begin{proof}[Proof of Theorem \ref{ogt2}]
It is sufficient to show that 
(\ref{ogt1h2}) holds true and to use
Theorem \ref{ogt1}. We assume that
$(u_1-u_2)^-$ belongs to $L^\infty(\OU)$. 
\par
According to the properties of the difference of two renormalized
solutions (see \cite{DMOP99}) we have for all $K>0$
\begin{gather}\label{ogt2n1}
  \int_{\{|u_1-u_2|<K\}}
  \big(\aop(x,Du_1)-\aop(x,Du_2)\big)\cdot(Du_1-Du_2) \dx \leq C K,
\end{gather}
  where $C$ is a constant independent of $K$. 
\par
Let $M$ be a real number such that
$M>\|(u_1-u_2)^-\|_{L^\infty(\OU)}$ and let $K>0$, $n\in\N^*$ and $\OU_{n,K}$
the set defined by (\ref{ogt1n11}). Since $\OU\subset \{ -M<u_1-u_2\}$
we get $\OU_{n,K}\subset \{|u_1|<2n\}\cap\{|u_2|<2n\}\cap \{
|u_1-u_2|<\max(M,K)\} $ and therefore
\begin{equation*}
  \frac{1}{n} \int_{\OU_{n,K}} |Du_1-Du_2|^p \dx 
\leq \frac{1}{n} \int_{ \substack{\{|u_1|<n,\, |u_2|<n\} \\ \cap {
      \{|u_1-u_2|<\max(M,K)\}}}}  |Du_1-Du_2|^p \dx.  
\end{equation*}
In both cases ($p<2$ and $p\geq 2$), the strong monotonicity of the
operator $\aop$, H\"older inequalities together with (\ref{ogrn1})
(with $\varphi\equiv 1$) and (\ref{ogt2n1}) 
 allow us to prove that for all $K>0$
\begin{gather*}
  \lim_{n\rightarrow+\infty } \frac{1}{n} \int_{ \substack{\{|u_1|<n,\,
      |u_2|<n\} \\ \cap {
      \{|u_1-u_2|<\max(M,K)\}}}}  |Du_1-Du_2|^p \dx = 0.
\end{gather*}
It follows that the conditions of Theorem \ref{ogt1} are satisfied and
then $u_1=u_2$.
\end{proof}
\par

\begin{rmk}
  In Theorem \ref{ogt1}, assuming $\aop$ to be strongly monotone,
  if condition (\ref{ogt1h2}) is satisfied for $K=0$ only (and not for
  every $K>0$), then
  $u_1=u_2$. Indeed, in this case (\ref{ogrn1}), (\ref{ogt2h}) and
  (\ref{ogt2n1}) imply that
$    \lim_{n\rightarrow+\infty}\frac{1}{n}\int_{\{|u_1-u_2|<K\}}
  |Du_1-Du_2|^p\dx=0$ $\forall K>0$ and then (\ref{ogt1h2}) is satisfied for all $K>0$.
 \end{rmk}
\par
\medskip
{\noindent \slshape Acknowledgments: } The author thanks F. Murat
  for interesting discussions and remarks regarding this paper.
\par

\bibliographystyle{plain}

\end{document}